\documentclass{article}
\usepackage{amssymb}
\title{{\bf From the representation theory of vertex operator algebras
to modular tensor categories in conformal field theory}}
\author{James Lepowsky}
\date{(Commentary on Yi-Zhi Huang's companion PNAS article ``Vertex operator
algebras, the Verlinde conjecture, and modular tensor categories'')}
\begin{document}
\bibliographystyle{plain}
\maketitle

Two-dimensional conformal quantum field theory (CFT) has inspired an
immense amount of mathematics and has interacted with mathematics in
very rich ways, in great part through the mathematically dynamic world
of string theory.  One notable example of this interaction is provided
by Verlinde's conjecture: E. Verlinde (1) conjectured that certain
matrices formed by numbers called the ``fusion rules'' in a
``rational'' CFT are diagonalized by the matrix given by a certain
natural action of a fundamental modular transformation (essentially, a
certain distinguished element of the group of two-by-two matrices of
determinant one with integer entries).  His conjecture led him to the
``Verlinde formula'' for the fusion rules and, more generally, for the
dimensions of spaces of ``conformal blocks'' on Riemann surfaces of
arbitrary genera.  A great deal of progress has been achieved in
interpreting and proving Verlinde's (physical) conjecture and the
Verlinde formula in mathematical settings, in the case of the
Wess-Zumino-Novikov-Witten models in CFT, which are based on affine
Lie algebras.  On the other hand, G. Moore and N. Seiberg (2, 3)
showed, on a physical level of rigor, that the general form of the
Verlinde conjecture is a consequence of the axioms for rational CFTs,
thereby providing a conceptual understanding of the conjecture.  In
the process, they formulated a CFT analogue, later termed ``modular
tensor category'' (cf. 4, 5) by I. Frenkel, of the classical notion of
tensor category for representations of (modules for) a group or a Lie
algebra.  It remained a very deep problem to construct, in a
mathematical as opposed to physical sense, structures (``theories'')
satisfying these axioms for rational CFT.  These axioms are, in fact,
much stronger than the Verlinde conjecture and modular tensor category
structure, and indeed, the mathematical construction of CFTs (as
opposed to the physical assumption that they should exist) is a very
rich field of study to which many mathematicians have contributed.  In
this issue of PNAS, Y.-Z. Huang (6) announces a (mathematical) proof
of the Verlinde conjecture in a very general form, along with two
notable consequences: the rigidity and modularity of a previously
constructed tensor category.  This work, along with results used in
this work, includes in particular the (mathematical) construction of a
significant portion of CFT---structures that actually do satisfy the
axioms---using the representation theory of vertex operator algebras.

Huang (6) invokes a great deal of earlier work in vertex (operator)
algebra theory.  The mathematical foundation of CFT may be viewed as
resting on the theory of vertex operator algebras ((7); see also (8)),
which reflect the physical features codified by
Belavin--Polyakov--Zamolodchikov (9).  Mathematically, vertex operator
algebra theory is extremely rich.  For the work discussed here, one
needs the representation theory of vertex operator algebras,
especially a tensor product theory for modules for a suitable vertex
operator algebra.  In classical tensor product theories for modules
for a group or for a suitable algebra such as a Lie algebra, one
automatically has the tensor product vector space available, and one
endows it with tensor product module structure by means of a natural
coproduct operation.  A module map from the tensor product of two
modules to a third module then amounts to an ``intertwining operator''
satisfying a natural condition coming from the group or algebra
actions on the three modules.

However, vertex operator algebra theory is imbued with considerable
``non-classical'' subtleties, intimately related in fact to the
non-classical nature of string theory in physics, and to construct a
tensor product theory of modules for a vertex operator algebra, one is
forced to proceed ``backwards'': First, one defines suitable
``intertwining operators'' (3, 10) among triples of modules.  The
dimensions of the spaces of these intertwining operators are the
fusion rules referred to above.  Then one has to construct a tensor
product theory that ``implements'' these intertwining operators.
After work of Kazhdan--Lusztig (11) for certain structures based on
affine Lie algebras, a tensor product theory for modules for a
suitable, general vertex operator algebra was constructed in a series
of papers summarized in (12).  Elaborate use of ``formal calculus''
(cf. (8)) was required in this work.  The main paper in this series
(13) establishes Huang's associativity theorem, which leads quickly to
the (categorical) coherence of the resulting braided tensor category.
In CFT terminology, this associativity theorem asserts the existence
and associativity of the operator product expansion for intertwining
operators, an assertion that was a key assumption (not theorem) in
(3).  In addition to braided tensor category structure, this series of
papers constructs the much richer ``vertex tensor category''
structure, which involves the conformal-geometric structure
established in (14), on the module category of a suitable vertex
operator algebra.

A fundamental theorem establishing natural modular transformation
properties of ``characters'' of modules for a suitable vertex operator
algebra was proved by Y. Zhu (15).  Requiring all of the theory
mentioned here, as well as results in (16)--(18), Huang formulates a
general, mathematically precise, statement of the Verlinde conjecture
in the framework of the theory of vertex operator algebras.  Assuming
only such purely algebraic, natural hypotheses as simplicity of the
vertex operator algebra, complete reducibility of suitable modules,
natural grading restrictions and cofiniteness, hypotheses that are
relatively easily checked and have indeed been previously verified in
a wide range of important families of examples, Huang skeches his
proof (see (6)).  The proof is heavily based on the results of his
recent papers (19) and (20), in which natural duality and modular
invariance properties for genus-zero and genus-one multipoint
correlation functions constructed from intertwining operators for a
vertex operator algebra satisfying the general hypotheses are
established; the multiple-valuedness of the multipoint correlation
functions leads to considerable subtleties that had to be handled
analytically and geometrically, rather than just algebraically.  The
strategy of the proof reflects the pattern of (2) and (3); in fact,
the main work is to establish two formulas of Moore and Seiberg that
they had derived from strong assumptions: the axioms for rational CFT.
The difficulties lie in the sequence of mathematical developments
briefly mentioned here.

As has been the case with many other major developments in the
mathematical study of string theory and conformal field theory over
the years, it is to be expected that the methods used by Huang (6)
will have further consequences.  In fact, this tensor category theory
has already been applied to a variety of fields in mathematics and
physics, including string theory or M-theory, in particular, D-branes.
The insight that continues to flow from the combined and respective
efforts of many physicists and mathematicians in this remarkable age
of string theory and its mathematical counterparts will surely produce
new surprises.

\bigskip

\noindent {\small \sc Department of Mathematics, Rutgers University,
Piscataway, NJ 08854}

\noindent {\em E-mail address}: lepowsky@math.rutgers.edu

\end{document}